\documentclass[12pt]{amsart}
 
\input{xy} 
\xyoption{all}

\newtheorem{theorem}{Theorem} [section]
\newtheorem{prop}[theorem]{Proposition} 
\newtheorem{lemma}[theorem]{Lemma}

\numberwithin{equation}{section} 
\usepackage{graphicx}

\newcommand\C{{\mathbb C}}

\newcommand\N{ {\mathbb N}}

\newcommand\cM{\mathcal{M}}

\newcommand{\cP}{\mathcal{P}}
\newcommand{\cD}{\mathcal{D}}

\renewcommand\phi{\varphi}

\newcommand\iso{\simeq} 

\renewcommand\mod{\operatorname{mod}}  

\newcommand{\cS}{\mathcal{S}}
\newcommand\ord {\operatorname{ord}} 
\newcommand\TF {\operatorname{TF}}
\newcommand\SF {\operatorname{SF}}
\newcommand\Spines {\operatorname{Spines}}
\newcommand\Top {\operatorname{Top}}

\begin{document}

\title[The geometry of the critically-periodic curves]{The geometry of the critically-periodic curves \\ in the space of cubic polynomials}

\author{Laura De Marco}
\address{Department of Mathematics, Statistics, and Computer Science,  University of Illinois at Chicago.}
\email{demarco@math.uic.edu}

\author{Aaron Schiff}

\begin{abstract}  
We provide an algorithm for computing the Euler characteristic of the curves $\cS_p$ in $\cP_3^{cm}\iso \C^2$, consisting of all polynomials with a periodic critical point of period $p$ in the space of critically-marked, complex, cubic polynomials.  The curves were introduced in \cite{Milnor:cubics, Bonifant:Kiwi:Milnor}, and the algorithm applies the main results of \cite{DP:combinatorics}.  The output is shown for periods $p \leq 26$.  
\end{abstract}

\date{\today}
\thanks{Research of both authors supported by the National Science Foundation.}

\maketitle

\thispagestyle{empty}

\section{Introduction}

Let $\cP_3^{cm}$ denote the space of cubic polynomials with marked critical points.  It is convenient to parametrize the space $\cP_3^{cm}$ by $(a,v)\in \C^2$, where the pair $(a,v)$ corresponds to the polynomial 
	$$f_{a,v}(z) = z^3 - 3a^2 z + 2 a^3 + v$$
with critical points at $\pm a$ and critical value $v = f_{a,v}(+ a)$.  

In this article, we study the geometry of the curves $\cS_p\subset\cP_3^{cm}$, introduced by J.~Milnor in \cite{Milnor:cubics},  consisting of cubic polynomials $f_{a,v}$ for which the critical point $+ a$ has period exactly $p$.  That is, 
	$$\cS_p = \{ (a,v) \in \cP_d^{cm}: f_{a,v}^p(a) = a, \; f_{a,v}^k(a) \not=a \mbox{ for all } 1 \leq k < p \}.$$
The curve $\cS_p$ is smooth for all $p$ \cite[Theorem 5.1]{Milnor:cubics}.  As a (possibly disconnected) Riemann surface, the curve $\cS_p$ has finite type:  it is obtained from a compact Riemann surface $\overline{\cS_p}$ by removing finitely many points.  The punctures lie at infinity in the space $\cP_d^{cm}$.  To date, the irreducibility of $\cS_p$ is unknown, though it is shown in \cite[\S8]{Bonifant:Kiwi:Milnor} for periods $p \leq 4$.  

The goal of this article is to explain an algorithm to compute the Euler characteristic of the compactification $\overline{\cS_p}$.  In \cite[Theorem 7.2]{Bonifant:Kiwi:Milnor}, it is shown to satisfy
\begin{equation} \label{BKM}
	\chi(\overline{\cS_p}) =  d_p \, (2-p) + N_p.
\end{equation}
The number $d_p$ is the degree of the curve $\cS_p$, and it is easily computable from the defining equation.  The number $N_p$ denotes the number of ends of $\cS_p$, the punctures $\overline{\cS_p}\setminus \cS_p$.  Our contribution is the algorithmic process to compute $N_p$, applying the main results of \cite{DP:combinatorics}.  The Euler characteristic of $\overline{\cS_p}$ is shown in Table \ref{output}, to period $p=26$.

We remark that the computation of the Euler characteristic $\chi(\overline{\cS_p})$ cannot be handled by traditional methods beyond the small periods.  A quick genus computation with Maple{{\scriptsize\textsuperscript{TM}}}, for example, yielded Euler characteristics for $p\leq 4$ and failed to provide an output for $p=5$ where $\cS_5$ is a curve of degree 80.    The degree of $\cS_p$ is on the order of $3^{p-1}$, and the curves $\overline{\cS_p}$ will be highly singular at infinity for any choice of projective compactification of $\cP_3^{cm} \iso \C^2$ and $p$ sufficiently large.  The Euler characteristics for periods $p\leq 4$ appear in \cite{Bonifant:Kiwi:Milnor}.  

\begin{table}[h]
\begin{center}  \begin{small}
\begin{tabular}{ | c | c | c | c | c |}
	\hline & & & & \\
Period & Tau-functions & Central ends  & Euler characteristic & $-\chi(\overline{\cS_p})/3^{p-1}$  \\  

\hline
1	&	1	&	1	&	2	&	-2.000	\\
2	&	1	&	1	&	2	&	-0.667	\\
3	&	3	&	5	&	0	&	0.000	\\
4	&	6	&	13	&	-28	&	1.037	\\
5	&	15	&	41	&	-184	&	2.272	\\
6	&	29	&	109	&	-784	&	3.226	\\
7	&	69	&	341	&	-3236	&	4.439	\\
8	&	141	&	973	&	-11848	&	5.417	\\
9	&	308	&	2853	&	-42744	&	6.515	\\
10	&	649	&	8301	&	-147948	&	7.517	\\
11	&	1406	& 24533		&	-505876	&	8.560  \\
12	&	2969	& 71737	&	-1694848	&	9.568 \\
13	&	6400	& 211653	&	-5630092	&	10.594 \\
14	&	13636	& 623485	&	-18491088	&	11.598 \\
15	&	29284	& 1842585	&	-60318292	&	12.611 \\
16	&	62746	& 5447957	&	-195372312	&	13.616 \\
17	&	134966	& 16134965	&	-629500300	&	14.624 \\
18	&	290089	& 47820749	&	-2018178784 &	15.628 \\
19	&	625298	& 141888285 &	-6443997868	& 16.633 \\
20	&	1348264	& 421295297 &	-20498523376 & 17.637 \\
21	&	2912779	& 1251903973 &	-64995935796 & 18.641 \\
22	&	6298309	& 3722380213 &	-205481381144 &  19.644 \\
23	&	13639477	& 11074683701 & -647923373764 &  20.647 \\
24	&	29567647	& 32965853477 &-2038171671252 &  21.650 \\
25	&	64181452	& 98175789309 &	-6397686770076 &  22.652 \\
26	&	139464021	& 292501047833	&	-20042379058084  & 23.655 \\
\hline
\end{tabular} \bigskip
\caption{The output of the Euler Characteristic algorithm.  From left to right:  the period $p$; the number of tau-functions with period $p$; the number of escape regions  of $\cS_p$ with the hybrid class of $z^2$ (see Theorem \ref{ends}); the Euler characteristic $\chi(\overline{\cS_p})$; and a comparison to $3^{p-1}$.}
\label{output} 
\end{small} \end{center}  
\end{table}


\subsection{Outline of the algorithm}
As described in \cite{Milnor:cubics}, the ends of $\cS_p$ correspond to the {\em escape regions} of $\cS_p$, the open subsets of $\cS_p$ consisting of polynomials with the critical point $-a$ tending to infinity under iteration.  The main ingredient in the computation of $N_p$ is the combinatorial analysis of polynomial dynamics on the basin of infinity, developed in \cite{Branner:Hubbard:2} and \cite{DP:combinatorics}.  Recall that the basin of infinity of a polynomial $f$ is the domain 
	$$X(f) = \{z\in\C: f^n(z) \to \infty\}.$$
From \cite{Branner:Hubbard:2}, we use the properties of the tableau (or equivalently, the Yoccoz tau-function) of a cubic polynomial; this combinatorial object encodes the first-return of a critical point to its ``critical nest."  From \cite{DP:combinatorics} we use the combinatorics of the pictograph, a more refined encoding of the first-return of a critical point to a ``decorated critical nest," allowing us to distinguish and count topological conjugacy classes.

The steps of the algorithm are:
\begin{enumerate}
\item	
Fix $p$.  For each $k$ dividing $p$, with $1 \leq k \leq p$, determine all admissible tau-functions with period $k$. 
\item
Count the number of topological conjugacy classes of basins of infinity $(f, X(f))$ associated to each tau-function.  
\item	
Compute the number of topological conjugacy classes of polynomials in $\cS_p$ with one escaping critical point:  each class is determined by the class of its basin of infinity (with a tau-function of period $k$) and a point in the Mandelbrot set associated to a period $p/k$ critical point.  
\item
Determine the number $N_p$ of escape regions in $\cS_p$:  there are either one or two ends in $\cS_p$ associated to each topological conjugacy class computed in the previous step, determined by the twist period of the tau-function.   
\item
Test the output against the degree of $\cS_p$:  $N_p$ is the total number of escape regions, while the degree of $\cS_p$ must equal the number of escape regions counted {\em with multiplicity}.  The multiplicity is computed from the tau-function.
\end{enumerate}

\noindent
Step (1) uses the tableau rules of \cite{Branner:Hubbard:2}, as corrected in \cite{Kiwi:cubics, DM:trees}; a translation into the language of the Yoccoz tau-functions was given in \cite{DS:count}.  The bulk of the computing time and memory usage goes into Step (1).  In \S\ref{tau}, we provide the theoretical results needed for the computation.  We include the theoretical results we used for improving the speed of the algorithm; we believe that some of these are interesting in their own right.

Step (2) was implemented already in \cite{DS:count}, applying the results of \cite{DP:combinatorics}.  Step (3) relies on the work of Branner and Hubbard in \cite{Branner:Hubbard:2} (see also \cite[Theorem 3.9]{Bonifant:Kiwi:Milnor}), to know that the conformal class of a cubic polynomial in an escape region depends only on the class of its basin and the class of its degree 2 polynomial-like restriction.  Steps (4) and (5) are explained in \S\ref{escape}, where we relate an escape region in $\cS_p$ to its quotient in the moduli space of cubic polynomials $\cM_3^{cm}$.  The multiplicity of an escape region is computed and depends only on the underlying tau-function.

\subsection{Details of the computation.}
An implementation of the algorithm was written with C++.  We compiled the output in Table \ref{output} to period $p = 26$.  The low periods are computed quickly, while the computation for period 26 took 9:13 hours (Intel Core 2 Quad @ 2.5 GHz on Windows 7 32-bit edition), executed on a single thread.  

\subsection{The growth rate of $\chi(\overline{\cS_p})$}
An easy computation shows that $-\chi(\overline{\cS_p})\to \infty$ as $p\to\infty$ \cite{Milnor:cubics}.  Using methods from pluripotential theory, Dujardin showed that 
	$$\frac{-\chi(\overline{\cS_p})}{3^p} \to \infty$$
as $p\to\infty$ \cite{Dujardin:cubics}.  After viewing the output of this algorithm, Milnor asked whether we have 
\begin{equation} \label{growth}
	\frac{\chi(\overline{\cS_p})}{3^{p-1}} = -p + O(1)
\end{equation}
as $p\to \infty$.  Or, equivalently by equation (\ref{BKM}), do we have 
	$$N_p = O(3^{p-1}) ?$$
We include the ratio $-\chi(\overline{\cS_p})/3^{p-1}$ in Table \ref{output}.

\subsection{Acknowledgments}  We would like to thank Jan Kiwi, Jack Milnor, and Kevin Pilgrim for helpful conversations and attention to the output.


\bigskip

\section{The $\tau$ functions}
\label{tau}

In this section, we define the Yoccoz tau-function of a cubic polynomial and explain Step 1 of the algorithm, the procedure to compute all periodic tau-functions of a given period $p$.  The main theoretical result is the following:

\begin{theorem} \label{tau length}
For each period $p \geq 1$, a tau-function has period $p$ if and only if 
	$$\tau(n) = n-p$$
for all $n \geq 2p-2$.
\end{theorem}

\noindent
We show that the bound $2p-2$ is optimal: for every $p\geq 3$, there exists a (unique) period $p$ tau-function with $\tau(2p-3) \not= p-3$.  See Lemma \ref{exception 2p-2}.

As described below, it is quite easy (from a theoretical point of view) to generate the periodic tau-functions, combining Theorem \ref{tau length} with  Theorem \ref{tau extension}.   A first approach might be to generate {\em all} admissible tau-functions of length $2p-2$ and test for equality $\tau(2p-2) = p-2$. As witnessed by the computations of \cite{DS:count}, however, the number of tau-functions grows exponentially with length, and only a small proportion are periodic.  For example, there are 649 tau-functions of period $p=10$, while there are 279,415 tau-functions of length $2p-2 = 18$.  Much of this section is devoted to the results we apply to reduce the computation time and memory usage.

\subsection{The tau-function of a polynomial.} \label{tau definition}
Fix a cubic polynomial $f$ with disconnected Julia set, and let 
	$$G_f(z) = \lim_{n\to\infty} \frac{1}{3^n} \log^+|f^n(z)|$$
be its escape rate.  Let $c_1$ and $c_2$ be the critical points of $f$, labeled so that $G_f(c_2)\leq G_f(c_1)$.   For each integer $n\geq 0$ such that $G_f(c_2) < G_f(c_1)/3^{n-1}$, we define the {\em critical puzzle piece} $P_n(f)$ as the connected component of $\{z: G_f(z) < G_f(c_1)/3^{n-1}\}$ containing $c_2$.  The puzzle piece $P_0(f)$ contains both critical points.  For positive integers $n$, we set	
	$$\tau(n) = \max\{j< n: f^{n-j}(c_2) \in P_j(f)\},$$
defining a function $\tau$ from $\{1, \ldots, N\}$ (or all of $\N$) to the non-negative integers.  The largest $N$ on which $\tau$ is defined is said to be the {\em length} of the tau-function.  In other words, $N$ is the greatest integer such that $G_f(c_2) < G_f(c_1)/3^{N-1}$.  If there is no maximal $N$, we say $\tau$ has length $\infty$.  

The {\em markers} of a tau-function with length $N$ are the integers 
	$$\{m\in \{1, \ldots, N-1\}: \tau(m+1) < \tau(m)+1\}.$$
The {\em marked levels} of $\tau$ are all integers in the forward orbit of a marker: 
	$$\{l\geq 0: l = \tau^n(m) \mbox{ for marker } m \mbox{ and } n>0\} \cup \{0\};$$
we say 0 is marked even if there are no markers.   The positive marked levels coincide with the lengths of the columns in the Branner-Hubbard tableau.  In terms of the polynomial $f$, a level $l$ is marked if the orbit of the critical point intersects $P_l(f) \setminus \overline{P_{l-1}(f)}$.  We say $l$ is {\em marked by k} if the $k$-th iterate $f^k(c_i)$ lies in $P_l(f) \setminus \overline{P_{l-1}(f)}$.

\subsection{Properties of tau-functions}  \label{properties}
Let $\N$ denote the positive integers $\{1, 2, 3, \ldots\}$.
For any positive integer $N$, a function 
	$$\tau: \{1, 2, 3, \ldots, N\} \to \N\cup\{0\}$$
or a function 
	$$\tau: \N \to \N\cup\{0\}$$
is said to be {\em admissible} if it satisfies the following properties (A)--(E):
\begin{itemize}
\item[(A)]	$\tau(1) = 0$
\item[(B)]	$\tau(n+1) \leq \tau(n)+1$
\end{itemize}
From (A) and (B), it follows that $\tau(n) < n$ for all $n\in\N$; consequently, there exists a unique integer $\ord(n)$ such that the iterate $\tau^{\ord(n)}(n) = 0$.  
\begin{itemize}
\item[(C)]	If $\tau(n+1) < \tau^k(n) + 1$ for some $0 < k < \ord(n)$, 
			then $\tau(n+1) \leq \tau^{k+1}(n) + 1$.
\item[(D)]	If $\tau(n+1) < \tau^k(n) + 1$ for some $0 < k < \ord(n)$, 
			and if $\tau(\tau^k(n) + 1) = \tau^{k+1}(n) + 1$, 
			then $\tau(n+1) < \tau^{k+1}(n) + 1$.
\item[(E)]	If $\ord(n)>1$ and $\ord(\tau^{\ord(n)-1}(n) + 1)=1$, then $\tau(n+1)\not=0$.
\end{itemize}

A tau-function is admissible if and only if it is the tau-function of a cubic polynomial \cite[Proposition 2.1]{DS:count}.  The proof is by induction on $N$, applying the rules for admissible tableaux in \cite{Branner:Hubbard:2}.  Property (E) is another formulation of the ``missing tableau rule" (M4) appearing in [Ki] and [DM].

Let $k$ be the number of markers which appear in the orbit 
	$$N\mapsto \tau(N) \mapsto \ldots \mapsto \tau^{\ord(N)}(N)=0,$$ 
and label these $k$ markers by $l_1', l_2', \ldots, l_k'$ so that
	$$N = l_0' > l_1' > l_2' > \cdots > l_k' > 0.$$
For each $0\leq i \leq k$, let $l_i = \tau(l_i')$ so that 
	$$\tau(N) = l_0 > l_1 > \cdots > l_k \geq 0.$$

Properties (A)--(E) imply the following:

\begin{theorem}  \cite[Theorem 2.2]{DS:count} \label{tau extension}
Given an admissible tau-function $\tau$ of length $N$, an extension to length $N+1$ is admissible if and only if 
	$$\tau(N+1) = l_i + 1 \mbox{ for some } 0 \leq i \leq k$$
or $\tau(N+1)=0$ if $l_k>0$ or $k=0$.  
\end{theorem}
 
\noindent
Note, in particular, that $\tau(N+1) = \tau(N) + 1$ is always an admissible extension to length $N+1$.

\subsection{Periodic tau-functions}
For cubic polynomials with exactly one critical point in the basin of infinity, the tau-function will have infinite length.  An admissible tau-function $\tau: \N\to \N\cup\{0\}$ is {\em periodic with period p} if there exists $N(\tau) \in \mathbb{N}$ such that 
	$$\tau(n) = n-p$$
for all $n\geq N(\tau)$.  Such tau-functions correspond to basins of infinity with a bounded critical orbit in a periodic component of the filled Julia set; $\tau$ has period $p$ if and only if the component has period exactly $p$.  For computational purposes, we need a bound on $N(\tau)$ depending only on the period $p$.  The bound $N(\tau) \leq 2p-2$ is granted by Theorem \ref{tau length}, which we prove below.

\begin{lemma}  \label{period p}
If $\tau$ has period $p$, then $\tau(n) \geq n-p$ for all $n$.  Further, if $\tau(n_0) = n_0-p$ for some $n_0$, then $\tau(n) = n-p$ for all $n\geq n_0$.  
\end{lemma}

\proof
This follows easily from property (B). 
\qed

\begin{lemma}  \label{p-1}
If $\tau$ has period $p$, then $l\leq p-1$ for all marked levels $l$.  
\end{lemma}

\proof
Let $f$ be any cubic polynomial with a given periodic tau-function.  Label the critical points of $f$ as in \S\ref{tau definition}.  Without loss of generality, we may assume the critical point $c_2$ is periodic with period exactly $p$. 

Suppose $l$ is marked by iterate $k$, and assume first that $\tau(l) = 0$.  From Lemma \ref{period p}, we have $l = l - \tau(l) \leq p$.  The first return of $P_l$ to the critical nest occurs with $f^l(P_l) = P_0$.  Because it maps with degree 2, the iterates $f^l(c_2)$ and $f^{k+l}(c_2)$ must lie in the two distinct components of $\{G_f < G_f(c_1)\}$ inside $P_0$.  By periodicity, then, we must have $l < p$.  

More generally, we have that the first return of $P_l$ to the critical nest is $f^{l-\tau(l)}(P_l) = P_{\tau(l)}$, and $l-\tau(l) \leq p$.  As above, because the first return is with degree 2, the images $f^{l-\tau(l)}c_2$ and $f^{k+ l-\tau(l)}(c_2)$ cannot lie in the same component of $\{G_f < G_f(c_2)/3^{\tau(l)}\}$ within $P_{\tau(l)}$, while $c_2$ and $f^k(c_2)$ do lie in the same component.  Therefore $l-\tau(l) < p$.  

In addition, we must have $l-\tau^2(p) \leq p$, as this is the first level where the forward orbit of $c_2$ and $f^k(c_2)$ might come together.   If $l$ is not a marker, then $f^{l-\tau(l)}(c_2)$ lies in the same component as $c_2$ at $\tau(l)$, and therefore its image at $\tau^2(l)$ is in a distinct component from that of $f^{k + l - \tau^2(l)}(c_2)$.  On the other hand, if $l$ is a marker, then $\tau(l)$ is marked by $l-\tau(l)$.  By periodicity, we can take $k = p- (l-\tau(l))$.  At $\tau^2(l)$, we have $f^{l-\tau^2(l)}(c_2)$ and $f^{k + l - \tau^2(l)}(c_2)$ again in distinct components.  In either case, we conclude that $l-\tau^2(l) < p$. 

We continue inductively.  For the induction step, we begin with $l - \tau^n(p) < p$ and $l-\tau^{n+1}(l) \leq p$.  We observe that at level $\tau^{n-1}(l)$, either $f^{l-\tau^{n-1}(l)}(c_2)$ or $f^{k+l-\tau^{n-1}(l)}(c_2)$ lies in the same component as $c_2$.  We consider the two cases:  if $\tau^{n-1}(l)$ is not a marker, then we may proceed two iterates to $\tau^{n+1}(l)$ keeping the image components distinct.  If $\tau^{n-1}(l)$ is a marker, then $\tau^{n-1}(l)$ is marked by $p - (\tau^{n-1}(l)-\tau^n(l))$; the component containing $f^{p - (\tau^{n-1}(l)-\tau^n(l))}(c_2)$ and the component containing $c_2$ at $\tau^{n-1}(l)$ must have distinct {\em preimages} at level $l$ which are sent to distinct components of $\tau^n(l)$, one of which contains $c_2$, and therefore to distinct components at $\tau^{n+1}(l)$.  We conclude that $l - \tau^{n+1}(l) < p$.  

Continuing until $\tau^{\ord(l)}(l) = 0$ completes the proof that $l < p$.  
\qed

\begin{lemma}  \label{p-2}
If $\tau$ has period $p$, and if a level $l$ is marked by $k = p-1$, then $l \leq p-2$.
\end{lemma}

\proof
Suppose $l$ is marked by $p-1$.  From Lemma \ref{p-1}, $l \leq p-1$.  By periodicity, $\tau(l) = l-1$.  From the admissible $\tau$ rules, it follows that $\tau(n) = n-1$ for all $1 \leq n \leq l$.  It follows that $n$ cannot be a marker for any $n \leq l-1$.  Consequently, level $n$ is marked by $l-n$ for all $0 \leq n \leq l-1$; in particular, $l$ marks level $0$.  Therefore $l \not= p-1$, because $p-1$ marks level $l$.


\medskip
\noindent{\bf Proof of Theorem \ref{tau length}.}
Suppose $\tau$ is periodic with period $p$.  By definition, there exists $N(\tau)$ so that $\tau(n) = n-p$ for all $n\geq N(\tau)$.  From Lemma \ref{p-1}, there are no marked levels $l\geq p$.  Therefore, there are no markers at levels $l \geq p + p-1 = 2p-1$.  Consequently, $\tau(n+1) = \tau(n) + 1$ for all $n\geq 2p-1$, and so we must have $\tau(n) = n-p$ for all $n\geq 2p-1$.  If $2p-2$ is a marker, then $\tau(2p-2) = p-1$, but this would imply that level $p-1$ is marked by $p-1$, violating Lemma \ref{p-2}.  Therefore, $\tau(2p-2) = p-2$.  
\qed

\begin{lemma}  \label{unique extension}
Suppose $\tau$ has length $N$ and $\tau(N) = N-p$.  Then $\tau$ extends uniquely to a sequence of period $p$, by setting 
	$$\tau(n) = n-p$$
for all $n> N$.  
\end{lemma}

\proof
The existence of the extension follows directly from Theorem \ref{tau extension}; the uniqueness from property (B).
\qed

\begin{lemma}  \label{marker < p}
Let $\tau$ have period $p$, and suppose $\tau(n_0) > n_0-p$ and $\tau(n_0+1) = n_0 + 1 -p$.  Then there exists a marker $m < p$ so that $\tau(m) = n_0-p$.  
\end{lemma}

\proof
By periodicity, there is some iterate $k$ so that $\tau^k(n_0) = n_0 - p$.  By assumption, $k>1$.  Let $m = \tau^{k-1}(n_0)$.  Because $\tau(n_0+1) = n_0 -p + 1$, we have that $n_0$ is a marker, so $m$ is marked.  By Lemma \ref{p-1}, then, $m < p$.  We need to show $m$ is also a marker.  Indeed, $\tau(m+1) = \tau(\tau^{k-1}(n_0) + 1) \not= n_0 - p + 1$ by property (D).   
\qed

\subsection{Examples/Exceptions}  \label{exceptions}
As demonstrated in Theorem \ref{tau length}, all periodic tau-functions of period $p$ must satisfy $\tau(n) = n-p$ for all $n\geq 2p-2$.  In fact, most periodic tau-functions of period $p$ also satisfy $\tau(n) = n-p$ for all $n\geq 2p-5$.  The following lemmas provide a complete list of the exceptions.  In the lemmas, we express the tau-function as a sequence of the form $\tau(1), \tau(2), \tau(3), \cdots$.  We remark that these lemmas are not used in the algorithm for the Euler characteristic computation, but we include them for completeness.

\begin{lemma} \label{exception 2p-2}  For each period $p\geq 3$,
there is a unique periodic tau-function with $\tau(n) = n-p$ for all $n\geq 2p-2$ and $\tau(2p-3) \not= p-3$.  It is given by 
\begin{itemize}
\item	$0, 1, 2,  \cdots, p-3, 0, 1, 2, \cdots, p-2, p-2, p-1, p, \cdots.$
\end{itemize}
\end{lemma}

\proof
By Lemma \ref{marker < p}, there is a marker $m < p$ with $\tau(m) = p-3$.  Thus $m$ can only be $p-2$ or $p-1$.  Consequently, the tau-function must begin with $0, 1, 2, \ldots, (p-3)$ or with $0, 0, 1, \ldots, (p-3)$.  In the first case, Theorem \ref{tau extension} implies that it can only be extended as $0, 1, 2, \ldots, (p-3), 0$ with $\tau(2p-3) = p-2$ and $\tau(2p-2) = p-2$.  In the case of $0, 0, 1, 2, \ldots, (p-3)$, if $p$ is even, then Theorem \ref{tau extension} implies the extension must be as $0, 0, 1, \ldots, (p-3), 1, 2, \ldots$, with $\tau(2p-3) = p-2$, but we cannot extend by $\tau(2p-2) = p-2$.  If $p$ is odd, then we must have $\tau(p) = 0$, but then $\tau(n) = n-p$ for all $n\geq p$.  
\qed

\begin{lemma} \label{exception 2p-3}  For each period $p\geq 4$, 
the only periodic tau-functions with $\tau(n) = n-p$ for all $n\geq 2p-3$ and $\tau(2p-4) \not= p-4$ are
\begin{itemize}
\item 	$0, 1, 2, \cdots, p-4, 0, 0, 1, 2, \cdots, p-3, p-3, p-2, p-1, \cdots$;  
\item		$0, 1, 2, \cdots, p-4, 0, 1, 2, \cdots, p-3, p-3, p-3, p-2, p-1, \cdots$
\end{itemize}
and if $p$ is odd then also
\begin{itemize}
\item		$0, 0, 1, 2, \cdots, p-4, 1, 2, \cdots, p-2, p-3, p-2, p-1, \cdots$.
\end{itemize}
\end{lemma}

\proof
By Lemma \ref{p-1}, we must have $\tau(2p-4) = p-3$ or $p-2$ or $p-1$.  Also, by Lemma \ref{marker < p}, level $p-4$ is marked by a marker $m< p$. 

Assume $\tau(2p-4) = p-3$.  Then $p-3$ is marked by $p-1$.  Periodicity implies that $p-4$ is marked by 1; that is $\tau(p-3) = p-4$.  The $\tau$ rules then imply that $\tau(n) = n-1$ for all $1 \leq n \leq p-3$, so our tau-function begins as $0, 1, 2, \ldots, (p-4)$.  Because $p-4$ must be marked, Theorem \ref{tau extension} implies that $\tau(p-2) = 0$.  Theorem \ref{tau extension} then allows for $\tau(p-1) = 0$ or $1$.  In either case, the tau-function is then uniquely determined by Theorem \ref{tau extension} and Lemma \ref{period p}, giving the first two possibilities stated in the Lemma. 

Now assume $\tau(2p-4) = p-2$.  Then level $p-2$ is marked by $p-2$, so by periodicity, we must have $\tau(p-2) = p-3$ or $\tau(p-2) = p-4$.  If $\tau(p-2) = p-3$, then the tau-function begins with $0, 1, 2, \ldots, p-3$, but then $p-4$ cannot be marked by a marker $m<p$ (contradicting Lemma \ref{marker < p}).  We must have $\tau(p-2) = p-4$, and the tau-function begins as $0, 0, 1, 2, \ldots, p-4$.   If $p$ is even, then we can only extend by 0 (for $p-4$ to be marked), but then $\tau(2p-4) \leq p-3$.  If $p$ is odd, then we can extend by $\tau(p-1) = 1$, and the final tau-function stated in the Lemma is admissible.  

The final possibility is that $\tau(2p-4) = p-1$.  The only way to mark $p-4$ by a marker $m<p$ is for $\tau(p-1) = p-4$, so the $\tau$ sequence begins with $0, \tau_2, \tau_3, 1, 2, \ldots, p-4$, for some $\tau_2, \tau_3 \leq 1$.  Then, as $p-4$ is marked by $p-1$, we must have $\tau(p) \leq 2$ by Theorem \ref{tau extension}.  But then $\tau(2p-4) < p-1$, so its $\tau$ orbit does not encounter any markers larger than 1, and we cannot have $\tau(2p-3) = p-3$.  
\qed

\begin{lemma}  \label{exception 2p-4} 
For each period $p\geq 5$, the only periodic tau-functions with $\tau(n) = n-p$ for all $n\geq 2p-4$ and $\tau(2p-5)\not=p-5$ are 
\begin{itemize}
\item		$0, 1, 2, \cdots, p-5, 0, 1, 2, \cdots, p-4, p-4, p-4, p-3, p-2, \cdots$;
\item		$0, 1, 2, \cdots, p-5, 0, 0, 1, 2, \cdots, p-4, p-4, p-4, p-3, p-2, \cdots$;
\item		$0, 1 ,2,\cdots,  p-5, 0, 0, 0, 1, 2, \cdots, p-4, p-4, p-3, p-2, \cdots$;
\item		$0, 1, 2 , \cdots, p-5, 0, 1, 0, 1, 2, \cdots, p-4, p-4, p-3, p-2, \cdots$;
\end{itemize}
and if $p$ is odd, then also
\begin{itemize}
\item		$0, 0, 1, 2, \cdots, p-5, 0, 1 , 2 ,\cdots, p-3, p-4, p-3, p-2, \cdots$; 
\end{itemize}
and if $p$ is even, then also
\begin{itemize}
\item		$0, 0, 1, 2, \cdots, p-5, 1, 1, 2, \cdots, p-3, p-4, p-3, p-2, \cdots$;
\end{itemize}
and if $(p-1)$ is divisible by 3, then also
\begin{itemize}
\item		$0, 1, 0, 1, 2, \cdots, p-5, 2, 3, \cdots, p-2, p-4, p-3, p-2, \cdots$; 
\end{itemize}
and if $(p-2)$ is divisible by 3, then also
\begin{itemize}
\item		$0, 0, 1, 1 ,2, \cdots, p-5, 2, 3, \cdots, p-2, p-4, p-3, p-2, \cdots$.
\end{itemize}
\end{lemma}

\proof
By Lemma \ref{marker < p}, $p-5$ is marked by a marker $< p$.  By Lemma \ref{p-1}, we must have $\tau(2p-5)$ equal to $p-4$, $p-3$, $p-2$, or $p-1$.  

Assume $\tau(2p-5) = p-4$.  Then level $p-4$ is marked by $p-1$.  Periodicity implies that $\tau(p-4) = p-3$, and therefore that $\tau(n) = n-1$ for all $n\leq p-4$.  Therefore, $\tau$ begins with $0, 1, 2, \ldots, p-5, 0$.  To reach $\tau(2p-5) = p-4$, we must have $1 \leq \tau(p) \leq 3$, allowing only the first four possibilities listed in the Lemma.  

Assume $\tau(2p-5) = p-3$.  We must have $\tau(p-3)$ equal to $p-4$ or $p-5$, by periodicity; for $p-5$ to be marked, we must have $\tau(p-3) = p-5$.  The tau-function begins with $0, 0, 1, 2, \ldots, p-5$.  If $p$ is odd, it can be continued by setting $\tau(p-2) = 0$, \ldots, $\tau(2p-5) = p-3$, and $\tau(2p-4) = p-4$.  If $p$ is even, then $\tau(p-2) = 1$, so we can take $\tau(p-1) = 1$ to allow for $\tau(2p-5) = p-3$.  

A similar argument handles the case of $\tau(2p-5) = p-2$.   
\qed

\bigskip
\section{Generating periodic tau-functions}

The goal is to generate a list of all tau-functions of period $p$.  For the later steps in the algorithm, we need the data of the tau-functions themselves, not only the total number.  

There is a unique tau-function of period $p=1$, given by $\tau(n) = n-1$ for all $n\geq 1$.  
For small periods, say period $p\leq 10$, there are few periodic tau-functions.  Applying Theorem \ref{tau extension}, we can generate all tau-functions to length $2p-2$.  From Theorem \ref{tau length}, the equality $\tau(2p-2) = p-2$ holds if and only if this tau-function extends to a sequence of period $p$; further, the extension is uniquely determined.  For example, the total number of tau-functions of length $8$ ($= 2p-2$ for $p=5$)  is only 144, so the computation time and memory usage are negligible for period $p=5$ \cite{DS:count}.  As the period grows, the total number of tau-functions of length $2p-2$ grows fast; it is probably larger than $4^{p-1}$.   We use the Lemmas of the previous section to reduce our computational requirements.

\subsection{Algorithm}

The algorithm proceeds as follows.  Fix $p> 1$.  We generate a list called {\em Periodic} containing all tau-functions of period $p$.  For the induction step, we generate a list called {\em Continue}.  

\medskip\begin{quote}
{\bf Initialization.}  Generate all tau-functions to length $n=p$, following Theorem \ref{tau extension}.  If $\tau(p) = 0$, include in {\em Periodic}.  If $\tau$ has no markers, then discard.  Otherwise, include in {\em Continue}.   
\end{quote}

\medskip\begin{quote}
{\bf Extension to length $n+1$ and test for periodicity.}   Choose $\tau$ from the list {\em Continue}.  Let $n$ be its length; by construction, $\tau(n) > n-p$.  Determine values $l_0, l_1, \ldots$ (as appearing in Theorem \ref{tau extension}, setting $N=n$) subject to the extra condition $n-p \leq l_i \leq l_0 = \tau(n)$.  For each such $l_i$, we consider the admissible extension of $\tau$, defined by 
	$$\tau(n+1) = l_i + 1.$$  
If $l_i  = n-p$, then include the extended $\tau$ in {\em Periodic}; by Lemma \ref{unique extension}, this $\tau$ uniquely determines a periodic tau-function.  

If $n < 2p-3$ and if $l_i > n-p$ and if 
	$$\max \{\tau(m):  m< p, \; \tau(m+1) \leq \tau(m)\} > n-p,$$
then include in {\em Continue}; this $\tau$ is a candidate to have a periodic extension, as it satisfies the necessary conditions of Lemmas \ref{period p} and \ref{marker < p} and Theorem \ref{tau length}.  Otherwise, discard.  Repeat the induction step until {\em Continue} is empty.  
\end{quote}

\subsection{Details}
In Tables \ref{table 10} and \ref{table 20}, we include the particulars of our computation for generating all periodic tau-functions of periods 10 and 20.  Following the algorithm above, we show the number of tau-functions in the lists {\em Periodic} and {\em Continue} as we increase the length of the tau-functions.  

\begin{table}[h]
\begin{center}  \begin{small}
\begin{tabular}{ | c | c | c | c | c |}
	\hline
Length & Periodic & Discard  & Continue \\  \hline
10 	&	205	&	1	&	435	  \\
11 	&  	201	&	242	&	506 		\\
12	&	139	&	567	&	479		\\
13	&	57	&	780	&	279		\\
14	&	26	&	497	&	134		\\
15	&	12	&	251	&	61		\\
16	&	6	&	122	&	21		\\
17	&	2	&	43	&	6		\\
18	&	1	&	13	&	0		\\
\hline
\end{tabular} \bigskip
\caption{Period 10 details:  generating the 649 tau-functions of period $10$ from a total of 279,415 tau-functions of length 18.  Final data file size = 7.8 KB, peak disk usage = 18 KB.}
\label{table 10} 
\end{small} \end{center}  
\end{table}

\begin{table}[h]
\begin{center}  \begin{small}
\begin{tabular}{ | c | c | c | c | c |}
	\hline
Length & Periodic & Discard  & Continue \\  \hline
20 	&	449308	&	1	&	848362		  \\
21	&	319756	&	528624	&	1055320		\\
22	&	389254	&	1059653	&	1116657		\\
23	&	114128	&	1523035	&	978211		\\
24	&	41925	&	1646071	&	674730		\\
25	&	17081	&	1299907	&	391444		\\
26	&	8896		&	800601	&	196937		\\
27	&	4138		&	403194	&	93346		\\
28	&	1898		&	192799	&	44601		\\
29	&	978	&	92478	&	20839		\\
30	&	475	&	43078	&	9636			\\
31	&	217	&	20028	&	4571			\\
32	&	113	&	9623		&	2054		\\
33	&	52	&	4309		&	932	\\
34	&	24	&	2004	&	414	\\
35	&	12	&	901	&	169	\\
36	&	6	&	373	&	57	\\
37	&	2	&	137	&	18	\\
38	&	1	&	41	&	0 \\
\hline
\end{tabular} \bigskip
\caption{Period 20 details:  generating the 1,348,264 tau-functions of period $p=20$ from a total of about 1.5 trillion tau-functions of length $2p-2 = 38$.  Final data file size = 29 MB, peak disk usage = 74 MB.}
\label{table 20} 
\end{small} \end{center}  
\end{table}

\bigskip
\section{Topological conjugacy classes of basins} \label{conjugacy classes}

In this section, we describe the algorithm to compute the number $\Top(\tau)$ of topological conjugacy classes of basins $(f, X(f))$ with a given tau-function $\tau$.  It is proved in \cite{DP:combinatorics} that $\Top(\tau)$ can be computed as
	$$\Top(\tau) = \Spines(\tau) \cdot \TF(\tau),$$
where $\Spines(\tau)$ is the number of {\em pictographs} (or {\em  truncated spines}) associated to $\tau$ and $\TF(\tau)$ is the associated {\em twist factor}.  We include here the steps to compute $\Spines(\tau)$ and $\TF(\tau)$.  These details already appeared in \cite{DS:count}.

The twist factor $\TF(\tau)$ is denoted by $\Top(\mathcal{D})$ in \cite{DP:combinatorics}, the number of conjugacy classes of basins with pictograph $\cD$, for any pictograph with tau-function $\tau$.  Indeed, it is easy to see that any pictograph with a period tau-function will have only finitely many marked levels, thus satisfying the hypotheses of \cite[Theorem 9.1]{DP:combinatorics}; further, it is stated there that the computation of $\Top(\cD)$ depends only on the underlying tau-function.

\subsection{Computing the number of pictographs}
Fix an admissible tau-function $\tau$ of length $N$.  As in \S\ref{tau definition}, the {\em markers} of $\tau$ are the integers 
	$$\{m\in \{1, \ldots, N-1\}: \tau(m+1) < \tau(m)+1\}.$$
The {\em marked levels} of $\tau$ are all integers in the forward orbits of the markers:
	$$\{l\geq 0: l = \tau^n(m) \mbox{ for marker } m \mbox{ and } n>0\} \cup \{0\};$$
we say 0 is marked even if there are no markers.    

As in Theorem \ref{tau extension}, we let $k$ be the number of markers which appear in the orbit 
	$$N\mapsto \tau(N) \mapsto \ldots \mapsto \tau^{\ord(N)}(N)=0.$$ 
Label these $k$ markers by $l_1', l_2', \ldots, l_k'$ so that
	$$N = l_0' > l_1' > l_2' > \cdots > l_k' > 0.$$
For each $0\leq i \leq k$, let $l_i = \tau(l_i')$ so that 
	$$\tau(N) = l_0 > l_1 > \cdots > l_k \geq 0.$$
For each $0 \leq i < k$, define $n_i$ by the condition that 
	$$\tau^{n_i} (l_i) = l_{i+1}$$
and define $n_k$ so that $\tau^{n_k}(l_k)=0$.  (The $n_i$ are called {\em special orders} in the program.)
For $0 < i < j \leq k+1$, we set 
	$$\delta(i,j) = \left\{ \begin{array}{ll} 
		1	&	\mbox{if } \tau(l_i'+1) = l_j+1 \\
		0	&	\mbox{otherwise} \end{array} \right. $$
where by convention we take $l_{k+1} = -1$.   Note that $\tau(l_k'+1) = 0$ for every $\tau$, so $\delta(k, k+1)=1$.  

The {\em symmetry} of $\tau$ is 
	$$s = \min\{n\geq 0: \tau^n(l_0) \mbox{ is a marked level}  \}.$$
Note that $s\leq n_0$.   To each admissible choice for $\tau(N+1)$ (from Theorem \ref{tau extension}) we define the $(N+1)$-th {\em spine factor} of $\tau$.  If $\tau(N+1) = l_i + 1$ with $i>0$, we set 
$$    \SF(\tau, N+1)  :=  
	2^{n_0-s} ( 2^{n_1} ( 2^{n_2} ( \cdots (2^{n_{i-1}} - \delta(i-1,i)) - \cdots ) - \delta(2, i) ) - \delta(1,i) );   $$
as above, we take $l_{k+1}=-1$.  If $\tau(N+1) = l_0+1 = \tau(N) + 1$, we set
	$$\SF(\tau, N+1) = 1$$
The number of pictographs (or equivalently, truncated spines) associated to a tau-function is computed inductively on the length.  

\begin{prop}  
Let $\tau$ be a periodic tau-function of period $p$.  The number of pictographs with tau-function $\tau$ is given by
	$$\Spines(\tau) = \prod_{j=1}^N \SF(\tau, j).$$
for any choice of $N$ with $\tau(N) = N-p$.
\end{prop}

\proof
That $\Spines(\tau)$ is the product of spine factors is deduced in \cite{DS:count}.  It remains to show that the computation terminates at a finite $N$ when $\tau$ is periodic.  From the definition of the spine factor, it is equal to 1 whenever $\tau(N+1) = \tau(N) + 1$.  For periodic taus, this will be the case for all $N$ sufficiently large.  Recall from Lemma \ref{period p} that once we find one $N$ with $\tau(N) = N-p$, this equality will hold for all $n\geq N$.  
\qed

\subsection{Computing the twist factor}  \label{twist factor}
Fix an admissible tau-function $\tau$ of length $N \in \N\cup\{\infty\}$ with finitely many marked levels.  For each $n < N$, the order of $n$ was defined in \S\ref{properties}; it satisfies $\tau^{\ord(n)}(n)=0$. For each marked level $l>0$, compute
	$$\mod(l) = \sum_{i=1}^l 2^{-\ord(i)}$$
and
	$$t(l) = \min\{n>0: n\mod(l) \in \N\}.$$
We define the {\em twist period} $T(\tau)$ by
\begin{equation} \label{T}
	T(\tau) = \max\{t(l): l \mbox{ is a marked level} \}
\end{equation}
or set $T(\tau) = 1$ if $\tau$ has no non-zero marked levels.  

Let $L(\tau)$ be the number of non-zero marked levels.  The {\em twist factor} is defined by
	$$\TF(\tau) = \frac{2^{L(\tau)}}{T(\tau)}.$$
Theorem 9.1 of \cite{DP:combinatorics} states that the number of topological conjugacy classes of basins associated to a given pictograph with tau-function $\tau$ is equal to $\TF(\tau)$.   

\subsection{The significance of the twist factor}  We include a few words here to explain the meaning of the values appearing in \S\ref{twist factor} to define the twist factor.  These play a role in the explanations of \S\ref{escape}.  

The quasiconformal deformations of the basin of infinity of a polynomial $f$ have a natural decomposition into twisting and stretching factors; see \cite{McS:QCIII} or the summary in \cite{DP:combinatorics}.  Let $f = f_{(a,v)}$ be a cubic polynomial with periodic tau-function $\tau$.  Let $G_f$ be its escape-rate function.  Recall that $-a$ is the critical point that escapes to infinity.  The {\em fundamental annulus} of $f$ is the domain 
	$$A(f) = \{z\in\C:  G_f(-a) < G_f(z) < 3G_f(-a)\}.$$
Viewing the basin of infinity $X(f)$ as an abstract Riemann surface, a full Dehn twist in $A(f)$ induces the {\em hemidromy} action described in \cite{Branner:Hubbard:2}; see also \cite{Branner:cubics} for an accessible summary.  

The twist period $T(\tau)$ is the least power of a full Dehn twist in the fundamental annulus that lies in the mapping class group of $f$.   To compute $T(\tau)$, we determine the induced amount of twisting in any image or preimage of $A(f)$ under the action of $f$.  Following the descriptions in \cite{Branner:cubics} and \cite{Branner:Hubbard:2}, it suffices to compute the relative moduli of these annuli lying between the two critical points;  the relative modulus of an annulus $A$ is the ratio $\mod(A)/\mod(A(f))$.  The value $\mod(l)$ computes exactly these sums of relative moduli down to the $l$-th marked level.  

The twist factor is the ingredient emphasized in \cite{DP:combinatorics}.  By measuring twist periods against the total number of ways to produce basins $(f, X(f))$ from a given pictograph, the discrepancy amounts to the twist factor $\TF(\tau)$.

\bigskip
\section{Escape regions}
\label{escape}

In this section we explain the final steps of the algorithm, incorporating the computations described in the previous section.

\subsection{The moduli space}
As discussed in \cite{Milnor:cubics}, there is a natural involution on the space $\cP_3^{cm}$, given by 
	$$I(a,v) = (-a,-v)$$
induced by the conjugation of $f_{(a,v)}$ by $z\mapsto -z$.  Thus there is a degree 2 projection
	$$\cP_3^{cm} \to   \cP_3^{cm}/I  =:  \cM_3^{cm}$$
to the moduli space of critically-marked cubic polynomials.  The action of $I$ preserves the curve $\cS_p$, defining a curve $\cS_p/I\subset \cM_3^{cm}$.

\subsection{Escape regions and multiplicity}
As introduced in \cite{Milnor:cubics} and \cite{Bonifant:Kiwi:Milnor}, an {\em escape region} of $\cS_p$ is a connected component of 
	$$\{ f_{(a,v)} \in \cS_p:  f^n_{(a,v)}(-a) \to \infty \mbox{ as }  n\to\infty \}.$$
That is, it consists of maps with one periodic critical point (at $+a$) and one escaping critical point (at $-a$).  

It follows from the general theory on stability that all polynomials in a given escape region $E$ are topologically conjugate on $\C$, as described in \cite{McS:QCIII}.  In this special setting, though, it can be seen directly from a canonical parameterization of $E$.  It is shown in \cite[Lemma 5.6]{Milnor:cubics} that each escape region $E$ is conformally a punctured disk, canonically identified with an $m$-fold cover of a punctured disk, for some positive integer $m = m(E)$.  This number $m(E)$ is called the {\em multiplicity} of $E$.

The covering map of degree $m(E)$ is defined by the assignment 
	$$(a,v) \mapsto \phi_{(a,v)}(2a),$$ 
where $\phi_{(a,v)}$ defines the uniformizing B\"ottcher coordinates near infinity for $f_{(a,v)}$, where $\phi_{(a,v)}(f_{(a,v)}(z)) = (\phi_{(a,v)}(z))^3$, unique if chosen to satisfy $\phi'(\infty)=1$.  The point $2a$ is the cocritical point for $-a$, so $f_{(a,v)}(2a) = f_{(a,v)}(-a)$.    In particular, the twisting deformation on the basin of infinity induces the change in angular coordinate on $E$.  In fact, the external angle of $2a$ is increased by $\pi$ under a full Dehn twist in the fundamental annulus of $f_{(a,v)}$; thus, $2 m(E)$ full twists closes a loop in $E$.

\begin{lemma}  \label{multiplicity}
Fix an escape region $E$ and let $\tau$ be the tau-function of any $f\in E$.  The multiplicity is given by 
	$$m(E) = \left\{  \begin{array}{ll} 1 & \mbox{if } T(\tau) = 1 \\
							T(\tau)/2 & \mbox{if } T(\tau) > 1 \end{array} \right.  $$
where $T(\tau)$ is the twist period computed in \S\ref{twist factor}.
\end{lemma}

\proof
Each escape region $E$ projects to an escape region $E/I$ in the curve $\cS_p/I \subset \cM_3^{cm}$.  By definition of the twist period, $T(\tau)$ full twists in a fundamental annulus are required to induce a closed loop in $E/I$.  But $E/I$ is doubly covered by a single escape region $E$ if and only if $f_{(a,v)}$ and $f_{(-a,-v)}$ are equivalent under a twist deformation, if and only if we have $T(\tau)=1$.  In this case of $T(\tau)=1$, two full twists are required to close a loop in $E$, corresponding to an argument increase of $2\pi$ for the cocritical point $2a$.  Therefore $m(E) = 1$.  On the other hand, if $T(\tau) > 1$, then each escape region $E$ projects bijectively to $E/I$;  thus $2 m(E) = T(\tau)$.  
\qed

\subsection{Hybrid classes}
For any polynomial $f$ in an escape region $E$ in $\cS_p$, the associated tau-function will have period $k$ for some $k$ dividing $p$.  A restriction of the iterate $f^k$ to a certain neighborhood of $+a$ will then define a quadratic {\em polynomial-like} map.  We refer to \cite{Douady:Hubbard:polynomial-like} for background information.  In this context, it is important to know that the conformal conjugacy class of $f$ is uniquely determined by the conformal conjugacy class of its basin of infinity $(f, X(f))$ and the {\em hybrid class} of its polynomial-like restriction \cite{Branner:Hubbard:2}.  See also \cite[Theorem 3.9, Corollary 3.10]{Bonifant:Kiwi:Milnor}.

We will use the following consequence of the general theory:  

\begin{prop}  \label{hybrid}
An escape region $E/I$ in $\cS_p/I$ is uniquely determined by 
\begin{enumerate}
\item		an integer $k$ dividing $p$ with $1 \leq k \leq p$; 
\item	  	a topological conjugacy class of basin dynamics $(f, X(f))$ with a critical end of period $k$; and 
\item		a point in the Mandelbrot set corresponding to a center of period exactly $p/k$.  
\end{enumerate}
\end{prop}

A center of period $n$ in the Mandelbrot set is a solution $c$ to the equation $f_c^n(0)=0$ where $f_c(z) = z^2+c$.  The center $c$ has period {\em exactly} $n$ if $n$ is the smallest positive integer for which the equality $f_c^n(0)=0$ holds. The number $\nu_2(n)$ of centers of period exactly $n$ is easily computable by the following relation:
		$$2^{n-1} = \sum_{q | n, \; 1 \leq q \leq n} \nu_2(q)$$

Combining the above results, we deduce the following:

\begin{theorem}  \label{ends}
For any tau-function $\tau$ with period $k$ dividing $p$, the number of escape regions in $\cS_p$ with tau-function $\tau$ is
	$$Ends(\tau, p) =  \left\{  \begin{array}{ll}  \nu_2(p/k) \, \Spines(\tau) \, \TF(\tau)  & \mbox{if } T(\tau) = 1 \\
							2\, \nu_2(p/k)\, \Spines(\tau) \, \TF(\tau) & \mbox{if } T(\tau) > 1 \end{array} \right.  $$
where $T(\tau)$ is the twist period, $\Spines(\tau)$ is the number of pictographs, and $\TF(\tau)$ is the twist factor of $\tau$.  The total number of escape regions in $\cS_p$ is therefore
	$$N_p = \sum_{k|p} \sum_{per(\tau) = k} \; Ends(\tau, p)$$
\end{theorem}

\noindent
In particular, in the case of $k=p$, $Ends(\tau, p)$ is the number of ``central ends" of $\tau$, coinciding with the number of all escape regions of $\cS_p$ with tau-function $\tau$ and hybrid class $z^2$.  The sum of $Ends(\tau, p)$ over all taus with period $p$ is shown in Table \ref{output}.  The sum of $Ends(\tau, p)$ over all taus with period {\em dividing} $p$ is the total number $N_p$ of escape regions in $\cS_p$.  

\proof
Fix $\tau$ of period $k$ dividing $p$.  From the arguments of \S\ref{conjugacy classes}, there are $\Spines(\tau) \, \TF(\tau)$ topological conjugacy classes of basins $(f, X(f))$ of cubic polynomials with tau-function $\tau$.  Applying Proposition \ref{hybrid}, there are consequently $ \nu_2(p/k)\, \Spines(\tau) \, \TF(\tau)$ escape regions $E/I$ in $\cS_p/I \subset \cM_3^{cm}$.  If $T(\tau) = 1$, then exactly as in the proof of Lemma \ref{multiplicity}, there is a unique escape region $E$ in $\cS_p$ mapped to each $E/I$.  If $T(\tau) >1$, there are exactly two escape regions mapped to each $E/I$.   
\qed

\subsection{Testing the computation}  \label{test}
We conclude with an explanation of the test of our computation against the degree of $\cS_p$.  

The multiplicity of an escape region $E$ in $\cS_p$ coincides with the number of intersection points of $E$ with any line in $\cP_3^{cm}$ of the form $\{a = a_0\}$ for any $a_0$ of sufficiently large modulus.  Therefore, the degree $d_p$ of the curve $\cS_p$ must satisfy 
	$$d_p = \sum_E m(E),$$
summing over all escape regions $E$ of $\cS_p$.   The degree $d_p$ is easily computed, as it satisfies:
	$$3^{p-1} = \sum_{q | p} d_q$$
where the sum is taken over all $q$ dividing $p$ with $1 \leq q \leq p$.  As established by Lemma \ref{multiplicity}, the value $m(E)$ depends only on the tau-function for the escape region $E$, so we may define 
	$$m(\tau) := m(E)$$
for any escape region $E$ associated to tau-function $\tau$.  

Our algorithm determines the value $Ends(\tau, p)$ for every tau-function of period $k$ dividing $p$; the ingredients are listed in Theorem \ref{ends}.  We can therefore check our computation by assuring equality of 
$$\sum_\tau \, m(\tau) \, Ends(\tau, p)  = d_p,$$
summing over all tau-functions $\tau$ of periods dividing $p$.

\bigskip\bigskip
\def\cprime{$'$}

 \end{document}